\pdfoutput=1 
\documentclass[draft=false, leqno]{article}
\usepackage{amsmath,amsthm,mathrsfs}
\usepackage[normal]{boilerart1}
\usepackage{longtable}

\usepackage[style=british]{csquotes}
\newcommand{\defn}[1]{\emph{#1}}

\newcommand{\comma}{\rlap{\ ,}}

\theoremstyle{definition}
\newtheorem*{ack*}{Acknowledgments}

\newcommand{\flowerstar}{f_{\ast}}

\newcommand{\cart}{\ensuremath{\textit{cart}}}

\newcommand{\coh}{\ensuremath{\textit{coh}}}

\newcommand{\cts}{\ensuremath{\textit{cts}}}

\newcommand{\equivalent}{\simeq}

\def\reflectedop#1#2{\mathop{\reflectbox{$#1{#2}$}}}

\newcommand{\Sh}{\categ{Sh}}
\newcommand{\PSh}{\categ{PSh}}
\newcommand{\Constr}{\categ{Constr}}

\newcommand{\Spacefin}{\Space_{\pi}}

\newcommand{\Lay}{\categ{Lay}}

\newcommand{\Aff}{\categ{Aff}}

\DeclareMathOperator{\Gal}{Gal}

\DeclareMathOperator{\Pro}{Pro}

\usepackage{todonotes}

\title{Exodromy for stacks}

\author{Clark Barwick \and Peter Haine}

\date{}

\begin{document}

\maketitle

\begin{abstract}
	In this short note we extend the Exodromy Theorem of \cite{exodromy} to a large class of stacks and higher stacks.
	We accomplish this by extending the Galois category construction to simplicial schemes.
	We also deduce that the nerve of the Galois category of a simplicial scheme is equivalent to its étale topological type in the sense of Friedlander.
\end{abstract}

\setcounter{tocdepth}{1}
\tableofcontents


\setcounter{section}{-1}
\section{Introduction}
In \cite{exodromy}, we identified a profinite category $\Gal(X)$ attached to any scheme\footnote{All our schemes and stacks in this paper will be assumed to be coherent.} $X$ \cites{Barwick:galperf}[Construction 13.5]{exodromy}.
The profinite category $ \Gal(X) $ classifies nonabelian constructible shea\-ves on $X$ (our \emph{Exodromy Equivalence} \cite[Theorem 11.7]{exodromy}) and the protruncated classifying space of $ \Gal(X) $ recovers the étale topological type of $X$ in the sense of Friedlander \cite{Haine:recon}.
A natural question, then, arises: what is the analogue of this construction for a simplicial scheme or stack?
For example, what is the correct exodromy representation corresponding to an equivariant constructible sheaf on a scheme with an action of a group scheme?

Here, we answer this question by extending the Galois category construction and the Exodromy Theorem to a large class of stacks and higher stacks.
Here is the basic construction.

\begin{cnstr}\label{cnstr:GalDelta}
	Let $Y_{\ast}$ be a simplicial scheme.
	Denote by $\Gal^{\Delta}(Y_{\ast})$ the following $1$-category.
	The objects are pairs $(m,\nu)$ consisting of an object $m \in \Delta$ and a geometric point $\nu \to Y_m$.
	A morphism $(m,\nu) \to (n,\xi)$ of $\Gal^{\Delta}(Y_{\ast})$ is a morphism $\sigma \colon m \to n$ of $\Delta$ and a specialisation $\nu \leftsquigarrow \sigma^{\ast}(\xi)$.
	This category has an obvious forgetful functor $\Gal^{\Delta}(Y_{\ast}) \to \Delta$, which is a cartesian fibration.
	A morphism $(m,\nu) \to (n,\xi)$ is cartesian over $\sigma \colon m \to n$ in $\Delta$ if and only if the specialisation $\nu \leftsquigarrow \sigma^{\ast}(\xi)$ is an isomorphism.

	The fibre over $m \in \Delta$ is the category $\Gal(Y_m)$, which we regard as a profinite category.
	(See \Cref{dfn:fibredinprofinite} for the precise notion of categories fibred in profinite categories.)
\end{cnstr}

Also attached to a simplicial scheme $Y_{\ast}$ is the étale topological type of $ Y_{\ast} $ as constructed by Eric Friedlander \cite[\S 4]{MR676809} and refined by David Cox \cite{MR550644}, Ilan Barnea and Tomer Schlank \cite{MR3459031}, David Carchedi \cite{Carchedi:higheretale}, and Chang-Yeon Cho \cite{MR3553672}.
The étale topological type of $ Y_{\ast} $ can be identified with the colimit in protruncated spaces of the simplicial object that carries $m \in \Delta$ to the protruncated étale homotopy type of $Y_m$ (see \cite[Theorem III.8]{MR550644}).
Since the protruncated homotopy type of the fibres of the cartesian fibration \smash{$\Gal^{\Delta}(Y_{\ast}) \to \Delta$} agree with the étale homotopy type of the schemes $Y_m$, it follows that the protruncated homotopy type of the the total category $\Gal^{\Delta}(Y_{\ast})$ is the colimit of this simplicial diagram.
In other words:

\begin{thm}
	The classifying protruncated space of $\Gal^{\Delta}(Y_{\ast})$ recovers the protruncated étale topological type of $Y_{\ast}$.
\end{thm}

This is a consequence of \Cref{prop:protruncated} below.
We will also show:

\begin{thm}[{\Cref{prop:GaldeltaConstr}}]
	If $Y_{\ast}$ is a presentation of an Artin $n$-stack $\mathscr{X}$, then the localisation of $\Gal^{\Delta}(Y_{\ast})$ at the cartesian edges classifies constructible sheaves on $\mathscr{X}$; in other words, a constructible sheaf on $\mathscr{X}$ is tantamount to a functor $\Gal^{\Delta}(Y_{\ast}) \to \Spacefin $ to $ \pi $-finite spaces that carries all cartesian edges to equivalences and restricts to a continuous functor $\Gal^{\Delta}(Y_m) \to \Spacefin$ for all $ m \in \Delta $.
\end{thm}

\noindent This theorem speaks only of Artin $n$-stacks, but it applies just as well to any coherent fpqc stack with a presentation as a simplicial scheme.

Additionally, this theorem speaks only about nonabelian constructible sheaves, but in fact the Galois categories we construct suffice to recover constructible $\overline{\QQ}_{\el}$ sheaves as well.
The proof will appear in a forthcoming note \cite{proexodromy}.

\begin{exm}
	Let $G$ be an affine group scheme over a ring $k$, and let $X$ be a $k$-scheme with an action of $G$.
	Then we have the usual simplicial $k$-scheme $B_{k,\ast}(X,G,k)$ whose $n$-simplices are $X \times_k G^{n}$;
	this presents the quotient stack $X/G$.

	Thus the category of $G$-equivariant (nonabelian) constructible sheaves on $X$ is equivalent to the category of continuous functors
	\[
		\Gal^{\Delta}(B_{k,\ast}(X,G,k)) \to \Spacefin
	\]
	that carry the cartesian edges to equivalences.
	If $\Lambda$ is a ring, then the derived category of $G$-equivariant constructible sheaves of $\Lambda$-modules on $X$ is equivalent to the category of continuous functors
	\[
		\Gal^{\Delta}(B_{k,\ast}(X,G,k)) \to \Perf(\Lambda)
	\]
	that carry cartesian edges to equivalences.

	The objects of the category $\Gal^{\Delta}(B_{k,\ast}(X,G,k))$ can be thought of as tuples
	\begin{equation*}
		(m, \Omega, x_0, g_1, \dots, g_m)
	\end{equation*}
	in which $m\in\Delta$ is an object, $\Omega$ is a separably closed field, and $x_0 \colon \Spec \Omega \to X$ and $g_1,\ldots,g_m \colon \Spec \Omega \to G$ are points with the property that $(x_0, g_1, \dots, g_m)$ is a geometric point of $X \times_{k} G^m$, so that $\Omega$ is the separable closure of the residue field of the image of the $(x_0, g_1, \dots, g_m)$ in the Zariski space of $ X \times_{k} G^m $.
\end{exm}

\begin{ack*}
	The second-named author gratefully acknowledges support from both the \textsc{mit} Dean of Science Fellowship and \textsc{nsf} Graduate Research Fellowship.
\end{ack*}


\section{Fibred {G}alois categories}

\begin{nul}\label{nul:conventionshighercats}
	We use the language and tools of higher category theory, particularly in the model of \emph{quasicategories}, as defined by Michael Boardman and Rainer Vogt and developed by André Joyal and Jacob Lurie.
	We will generally follow the terminological and notational conventions of Lurie's trilogy \cites{HTT,HA,SAG}, but we will simplify matters by \emph{systematically using words to mean their good homotopical counterparts.}
	So `category' here means `$\infty$-category', `topos' means `$\infty$-topos', \& c.

	We write $ \Space $ for the category of spaces and $ \Spacefin \subset \Space $ for the full subcategory spanned by the $ \pi $-finite spaces.

	We use \HTT{Corollary}{3.2.2.13} systematically to construct cartesian fibrations; we leave the details of this by now standard construction implicit in what follows.
\end{nul}

\begin{ntn}
	If $\XX \to S$ is a topos fibration \HTT{Definition}{6.3.1.6}, then for any morphism $f \colon s \to t$ of $S$, there is a corresponding geometric morphism $f_{\ast} \colon \XX_t \to \XX_s$ of topoi; its left exact left adjoint will be denoted $f^{\ast}$.
\end{ntn}

\begin{dfn}
	Let $S$ be a category.
	A \emph{bounded coherent topos fibration} $ \XX \to S $ is a topos fibration in which each fibre $\XX_s$ is bounded coherent, and for any morphism $f \colon t \to s$ of $S$, the induced geometric morphism $ \flowerstar \colon \XX_s \to \XX_t$ is coherent \cites[Definitions \SAGthmlink{A.2.0.12} \& \SAGthmlink{A.7.1.2}]{SAG}[Definition 5.28]{exodromy}.
	A \emph{spectral topos fibration} $ \XX \to S $ is a bounded coherent topos fibration in which each fibre $\XX_s$ is a spectral topos (for the canonical profinite stratification \cite[Lemma 9.40 \& Definition 10.3]{exodromy}).
\end{dfn}

\begin{nul}
	The usual straightening/unstraightening equivalence restricts to an equivalence between the category of bounded coherent (respectively, spectral) topos fibrations $\XX \to S$ and the category of functors from $S^{\op}$ to the category of bounded coherent (resp., spectral) topoi (cf. \HTT{Proposition}{6.3.1.7}).

	For a bounded coherent topos fibration $ \fromto{\XX}{S} $ we write $\XX^{\coh}_{<\infty} \subseteq \XX$ for the full subcategory spanned by the objects that are truncated and coherent in their fibre \cite[Definition 5.18]{exodromy}.
	Then $\XX^{\coh}_{<\infty} \to S$ is a cocartesian fibration that is classified by a functor from $ S $ to the category of bounded pretopoi \cite[\SAGthm{Definition}{A.7.4.1} \& \SAGthm{Theorem}{A.7.5.3}]{SAG}.
\end{nul}

\begin{exm}
	If $X_{\ast}$ is a simplicial (coherent!) scheme, then the fibred topos $X_{\ast,\et} \to \Delta$ is a spectral topos fibration.
\end{exm}

\begin{nul}\label{nul:Hochsterdual}
	Hochster duality \cite[Theorem 10.10]{exodromy} expresses an equivalence between the category of profinite layered categories\footnote{A category $ C $ is \defn{layered} if every endomorphism in $ C $ is an equivalence.} and the category of spectral topoi, which carries a profinite layered category $\Pi= \{\Pi_{\alpha}\}_{\alpha \in A}$ to the spectral topos $\widetilde{\Pi}$ of sheaves in the effective epimorphism topology \cite[\SAGsubsec{A.6.2}]{SAG} on the bounded pretopos
	\[
		\Fun^{\cts}(\Pi, \Spacefin)\coloneq\colim_{\alpha\in A^{\op}}\Fun(\Pi_{\alpha},\Spacefin)
	\]
	of \textit{continuous functors} $ \fromto{\Pi}{\Spacefin}$.
	Under Hochster duality, the category of spectral topos fibrations $\XX \to S$ is equivalent to the category of functors from $S^{\op}$ to the category of profinite layered categories.
\end{nul}

A fibred form of Hochster duality is what allows us to construct fibred Galois categories.
To define it, we need to make sense categories fibred in profinite stratified spaces.

\begin{dfn}\label{dfn:fibredinprofinite}
	Let $S$ be a category.
	A functor $\Pi \to S$ will be said to be a \emph{category over $S$ fibred in layered categories} if it is a catesian fibration whose fibres are layered categories.
	We write $\Lay^{\cart}_{/S}$ for the category of categories over $S$ fibred in layered categories.
\end{dfn}

\begin{cnstr}
	There is a monad $ T $ on the category $\Lay$ of small layered categories given by sending a layered category $\Pi$ to the limit over the $\pi$-finite layered categories to which it maps.\footnote{That is, $ T $ is the right Kan extension of the inclusion $ \incto{\Lay_{\pi}}{\Lay} $ of $ \pi $-finite layered categories along itself.}
	The category of $T$-algebras is equivalent to the category of profinite layered categories.
	If $S$ is a category, this monad can be applied fibrewise to give a monad $T_S$ on the category $\Lay^{\cart}_{/S}$ of categories fibred in layered categories.

	Under the straightening/unstraightening identification
	\begin{equation*}
		\Lay^{\cart}_{/S} \equivalent \Fun(S^{\op},\Lay) \comma
	\end{equation*}
	the monad $T_S$ corresponds to the monad on $\Fun(S^{\op},\Lay)$ given by applying $T$ objectwise.
	Consequently, the category of $T_S$-algebras is equivalent to the category of functors from $S^{\op}$ to the category of profinite layered categories.
\end{cnstr}

\begin{dfn}
	Let $S$ be a category.
	A \defn{category over $S$ fibred in profinite layered categories} is a $T_S$-algebra.
	If $\Pi \to S$ is a category fibred in layered categories, then a \defn{fibrewise profinite structure} on $\Pi \to S$ is a $T_S$-algebra structure on $\Pi \to S$.
	We write $\Lay^{\cart,\wedge}_{\pi,/S}$ for the category of $T_S$-algebras.
\end{dfn}

\begin{wrn}
	One might also contemplate the category $\Pro(\Lay^{\cart}_{\pi,/S})$ of proöbjects in the full subcategory
	\begin{equation*}
		\Lay^{\cart}_{\pi,/S} \subseteq \Lay^{\cart}_{/S}
	\end{equation*}
	spanned by those cartesian fibrations whose fibres are $\pi$-finite layered categories.
	This is generally \emph{not} equivalent to the category of categories over $S$ fibred in profinite layered categories.
	Under straightening/unstraightening, the category $\Lay^{\cart,\wedge}_{\pi,/S}$ is equivalent to the category $\Fun(S^{\op},\Lay^{\wedge}_{\pi})$, whereas $\Pro(\Lay^{\cart}_{\pi,/S})$ is equivalent to the category $\Pro(\Fun(S^{\op},\Lay_{\pi}))$.
	These coincide when $S$ is a finite poset \HTT{Proposition}{5.3.5.15}, but otherwise typically do not coincide.
\end{wrn}

\begin{nul}\label{nul:descriptionoffibredHochster}
	Let $S$ be a category.
	Then the category of spectral topos fibrations over $S$ is equivalent to the category $\Lay^{\cart,\wedge}_{\pi,/S}$.
	Let us make the equivalence explicit.
	If $\XX \to S$ is a spectral topos fibration, then we define a category over $S$ fibred in layered categories
	\begin{equation*}
		\Pi^{S,\wedge}_{(\infty,1)}(\XX) \to S
	\end{equation*}
	as follows.
	An object of $\Pi^{S,\wedge}_{(\infty,1)}(\XX)$ is a pair $(s,\nu)$, where $s\in S$ and $\nu_{\ast} \colon \Space \to \XX_s$ is a point.
	A morphism $(s,\nu) \to (t,\xi)$ is a morphism $f \colon s \to t$ of $S$ and a natural transformation $\nu_{\ast} \to f_{\ast}\xi_{\ast}$.
	The category $\Pi^{S,\wedge}_{(\infty,1)}(\XX)$ fibred in layered categories admits a canonical fibrewise profinite structure;
	the fibre $\Pi^{S,\wedge}_{(\infty,1)}(\XX)_s$ over an object $s\in S$ is the profinite stratified shape $\Pi^{\wedge}_{(\infty,1)}(\XX_s)$ of \cite[Construction 11.1]{exodromy}.

	In the other direction, if $\Pi \to S$ is a category over $S$ fibred in profinite layered categories, then let $X_0 \to S$ denote the cocartesian fibration in which the objects are pairs $(s,F)$ consisting of an object $s\in S$ and a functor $F \colon \Pi_s \to \Spacefin$, and a morphism $(f, \phi) \colon (s, F) \to (t, G)$ consists of a morphism $f \colon s \to t$ of $S$ and a natural transformation $\phi \colon f_!F \to G$.
	Then $(\widetilde{\Pi})^{\coh}_{<\infty}$ is equivalent to the subcategory of $X_0$ whose objects are those pairs $(s,F)$ in which $F$ is continuous and whose morphisms are those pairs $(f, \phi)$ in which $\phi$ is continuous \Cref{nul:Hochsterdual}.
\end{nul}

\begin{cnstr}
	If $S$ is a category and $\YY$ is a bounded coherent topos, then the projection $\YY \times S \to S$ is a bounded coherent topos fibration.
	The assignment $ \goesto{\YY}{\YY \times S} $ defines a functor from the category of bounded coherent topoi to the category of bounded coherent topos fibrations over $S$.
	This functor admits a left adjoint, which we denote by $|\cdot|_S$.
	At the level of pretopoi, $(|\XX|_S)^{\coh}_{<\infty}$ is equivalent to the category of cocartesian sections of $\XX^{\coh}_{<\infty} \to S$, i.e., the limit of the corresponding functor from $S$ to bounded pretopoi.
\end{cnstr}

Now we arrive at the main topos-theoretic result.
\begin{prp}\label{prop:maintopos}
	Let $S$ be a category, and let $\XX \to S$ be a spectral topos fibration.
	Then the pretopos $(|\XX|_S)^{\coh}_{<\infty}$ is equivalent to the category of functors $F \colon \Pi^{S,\wedge}_{(\infty,1)}(\XX) \to \Spacefin$ with the following properties.
	\begin{itemize}
	 	\item $F$ carries any cartesian edge to an equivalence.
	 	\item For any object $s\in S$, the restriction $F|_{\Pi^{\wedge}_{(\infty,1)}(\XX_s)}$ is continuous.
	 	\item $F$ is uniformly truncated in the sense that there exists an $N\in\NN$ such that for any object $(s,\nu) \in \Pi^{S,\wedge}_{(\infty,1)}(\XX)$, the space $F(s,\nu)$ is $N$-truncated.
	\end{itemize}
\end{prp}

\begin{proof}
	The pretopos $(|\XX|_S)^{\coh}_{<\infty}$ can be identified with the category of cocartesian sections of $\XX^{\coh}_{<\infty} \to S$.
	The description of \Cref{nul:descriptionoffibredHochster} completes the proof.
\end{proof}

Please note that the last condition of \Cref{prop:maintopos} is automatic if $S$ has only finitely many connected components (e.g., $ S = \Delta $).

\begin{exm}
	If $X_{\ast}$ is a simplicial scheme, then the category over $ \Delta $ fibred in profinite layered categories $ \Pi_{(\infty,1)}^{\Delta,\wedge}(X_{\ast,\et}) $ associated to the spectral topos fibration $X_{\ast,\et} \to \Delta$ is the category $\Gal^{\Delta}(X_{\ast})$ of \Cref{cnstr:GalDelta}.
	In this case, \Cref{prop:maintopos} implies that $(|X_{\ast,\et}|_{\Delta})^{\coh}_{<\infty}$ is equivalent to the category of functors $\Gal^{\Delta}(X_{\ast}) \to \Spacefin$ that carry cartesian edges to equivalences and restrict to continuous functors $\Gal^{\Delta}(X_m) \to \Spacefin$ for all $m \in \Delta$.
\end{exm}

Finally, since the profinite stratified shape is a delocalisation of the protruncated shape \cite[Theorem 2.5]{Haine:recon} we deduce the following:
\begin{prp}\label{prop:protruncated}
 	Let $S$ be a category, and let $\XX \to S$ be a spectral topos fibration.
 	Then the protruncated shape of $|\XX|_S$ is equivalent to the protruncated homotopy type of $\Pi^{S,\wedge}_{(\infty,1)}(\XX)$.
\end{prp} 

\begin{exm}
	If $X_{\ast}$ is a simplicial scheme, then the protruncated homotopy type of the fibrewise profinite category $\Gal^{\Delta}(X_{\ast})$ is equivalent to the Friedlander étale topological type of $X_{\ast}$ \cite[Theorem A]{Haine:recon}.
\end{exm}


\section{Sheaves on stacks}

\begin{cnstr}
	Write $ \Aff $ for the $ 1 $-category of affine schemes.
	We employ \HTT{Corollary}{3.2.2.13} to construct a category $ \PSh_{\et} $ and a cocartesian fibration
	\begin{equation*}
		\PSh_{\et} \to \Aff^{\op}
	\end{equation*}
	in which the objects of $\PSh_{\et}$ are pairs $(S, F)$ consisting of an affine scheme $S$ and a presheaf (of spaces) on the small étale site of $S$, and a morphism $(S, F) \to (T, G)$ is a pair $(f, \phi)$ consisting of a morphism $f \colon T \to S$ and a morphism of presheaves $\phi \colon f^{-1}F \to G$ on the small étale site of $T$.
	Define $\Sh_{\et} \subset \PSh_{\et}$ to be the full subcategory spanned by those pairs $(S, F)$ in which $F$ is a sheaf; then $\Sh_{\et} \to \Aff^{\op}$ is a topos fibration.
	Define $\Constr_{\et} \subset \Sh_{\et}$ to be the further full subcategory spanned by those pairs $(S, F)$ in which $F$ is a (nonabelian) constructible sheaf \cite[Definition 10.11]{exodromy}; then $ \Constr_{\et} \to \Aff^{\op}$ is a cocartesian fibration.
\end{cnstr}

\begin{dfn}
	Let $X \to \Aff$ be a stack, i.e., a right fibration that is classified by an accessible fpqc sheaf $\Aff^{\op} \to \Space$.
	A \emph{(nonabelian) constructible sheaf} on $X$ is a cocartesian section
	\[
		F \colon X^{\op} \to \Constr_{\et}
	\]
	over $\Aff^{\op}$.
	We write $ \Constr_{\et}(X) $ for the category of constructible sheaves on $ X $.
\end{dfn}

\begin{wrn}
	This can only be expected to be a reasonable definition for coherent stacks.	
\end{wrn}

\begin{nul}
	Informally, a constructible sheaf $F$ on $X$ assigns to every affine scheme $S$ over $X$ a constructible sheaf $F_S$ and to every morphism $f \colon S \to T$ of affine schemes an equivalence $F_S \simeq f^{\ast}F_T$.
	In other words, the category of constructible sheaves on $X$ is the limit of the diagram $X^{\op} \to \Cat$ given by the assignment $S \mapsto \Constr_{\et}(S)$.

	Of course, since $X$ is not a small category, it is not obvious that this limit exists in $\Cat$.
	However, if $X$ contains a small limit-cofinal full subcategory $Y$, then the desired limit exists.
\end{nul}

Now we conclude:

\begin{prp}\label{prop:GaldeltaConstr}
	If $p \colon X \to \Aff$ is a stack, and if $X$ is presented by a simplicial scheme $Y_{\ast}$, then we obtain an equivalence between the category $\Constr_{\et}(X)$ and the category of functors
	\[
		\Gal^{\Delta}(Y_{\ast}) \to \Spacefin
	\]
	that carry cartesian edges to equivalences and for all $m \in \Delta$ restrict to a continuous functor $\Gal(Y_m) \to \Spacefin$.
\end{prp}


\DeclareFieldFormat{labelnumberwidth}{#1}
\printbibliography[keyword=alph]
\addcontentsline{toc}{section}{References} 
\DeclareFieldFormat{labelnumberwidth}{{#1\adddot\midsentence}}
\printbibliography[heading=none, notkeyword=alph]

\end{document}